# On the properties of nodal price response matrix in electricity markets

Vadim Borokhov

*Abstract*-- We establish sufficient conditions for nodal price response matrix in electric power system to be symmetric and negative (semi-)definite. The results are applicable for electricity markets with nonlinear and intertemporal constraints.

*Index Terms*-- Power system economics, electricity market, nodal pricing, residual demand, price response.

## I. INTRODUCTION

Recent decades had witnessed gradual liberalization of electric power sectors in a number of markets with restructuring of vertically-integrated utilities into power generation, supply, and transmission/distribution businesses. Unbundling of competitive sectors from naturally monopolistic businesses paved the way to introduction of free market pricing for power. We will consider the power markets based on the nodal pricing [1]-[3] and focus on the sensitivity of the locational marginal prices (LMP) with respect to infinitesimal nodal injections of costless power. The sensitivity analysis of LMP naturally arises in a number of problems including profit maximization of a generation firm with market power [4].

Consider the firm $g$ operating generating units located in a set of nodes $\{i\}$. Let $x^{i,h}$ denote power injection of the corresponding generating unit in node (bus) $i$ at our $h$, $\lambda_{i,h}$ be the corresponding nodal price, $c(\{x^{i,h}\})$ be the firm power production cost function, then in the absence of binding constraints on variables $\{x^{i,h}\}$ the first-order conditions for profit optimization problem for $g$ selling all power at the DAM nodal prices implies

$$\left(\lambda_{i,h} - \partial c / \partial x^{i,h}\right) = -\sum_{i',h'} x^{i',h'} (\partial \lambda_{i',h'} / \partial x^{i,h}),$$

which gives

$$\sum_{i,h} \left(\lambda_{i,h} - \partial c / \partial x^{i,h}\right) x^{i,h} = -\sum_{i,i',h,h'} x^{i',h'} (\partial \lambda_{i',h'} / \partial x^{i,h}) x^{i,h}, \quad (1)$$

where the summation is performed over all nodes $i$, $i'$ of the network with $g$'s generating units and all hours $h$, $h'$ of the DAM period. Equation (1) implies that additional revenue received by the firm due to the exercise of the market power is set by the symmetric part of the nodal price response matrix $\partial \lambda_{i',h'} / \partial x^{i,h}$. If the matrix is negative semi-definite, then RHS of (1) is non-negative and the firm receives non-negative markup due to its market power. That is why symmetry and sign definiteness of the nodal price response matrix is of special interest for both market players and regulators.

The sensitivity problem for one-period DAM has been extensively studied in [5]-[9] for DC model, i.e., for the case of linear constraints without intertemporal constraints (such as ramping constraints). In [5]-[9] using sensitivity analysis of the first-order conditions of the corresponding optimization problem it was established that in DC model the nodal price response matrix is symmetric and negative (semi-)definite matrix. In [9] the closed form solution for the nodal price response matrix was obtained in DC model for a piece-wise quadratic cost function. The approach utilized in [5]-[9] heavily relies on the special structure of the first-order conditions in DC model and explicit expression for nodal price response matrix obtained from sensitivity analysis of the first-order conditions. For the case of linear constraints with quadratic power production cost function the price response matrix is constant in each interval with unaltered set of binding constraints.

Full optimal power flow model with set unit commitment schedule for multi-period DAM accounts for transit losses in electric lines, ramping constraints, etc., and hence includes nonlinear constraints as well as intertemporal constraints. In the case of nonlinear constraints, the price response matrix is not constant even for quadratic power production cost function (and may have significantly different sensitivity values for the large power systems with transit power losses), thus making the solution of the supplier profit maximization problem more challenging. Since binding intertemporal constraints may change the outcome of the nodal price response matrix calculation, they should be taken into account.

We study the multi-period DAM with nonlinear and intertemporal constraints (i.e., the cases when the explicit expression for the price response matrix is unavailable) and illustrate that symmetry and negative (semi-)definiteness of the matrix follow from more general conditions on the underlying optimization problem without use of explicit expression for the matrix. Thus, the main contribution of the present paper is to extend the results of [5]-[9] to the case of multi-period DAM with nonlinear and intertemporal constraints to show that the symmetry of nodal price response matrix originates from general properties of the Lagrange multipliers, when exogenous parameters enter additively in the binding constraints, satisfying the linear independence

---

Vadim Borokhov is with En+ Development, Shepkina 3, Moscow, 129090, Russia, (e-mail: vadimab@eurosib.ru). The views expressed in this paper are solely those of the author and not necessarily those of En+ Development.

constraint qualification (LICQ), while the negative (semi-)definiteness of the nodal price response matrix is readily obtained in this framework for the instance of convex constraints (which incorporates the DC model case as well). As in [9] we require LICQ because, in contrast to other more general constraint qualification conditions, it implies that the binding constraint set locally doesn't induce any constraints on the additive exogenous parameters, so that their values can vary independently and partial derivatives over the parameters can be considered.

The paper is organized as follows. Section II is devoted to application of the general statements to LMP, which are further generalized in section III, section IV provides illustration of the findings for the case with binding intertemporal constraint, the section V contains conclusion, in section VI we recall general properties of the Lagrange multipliers' derivatives with respect to parameters entering the constraints additively. Throughout the paper we adopt LICQ as constraint qualification condition.

## II. DAM WITH FIXED UNIT COMMITMENT SCHEDULE

Consider wholesale multi-period day ahead electric power market (DAM) with set unit commitment schedule operating on bid-based security constrained economic dispatch principle according to the financially binding offers/bids supplied by wholesale market players. Let DAM with hourly locational marginal pricing be cleared simultaneously for all hours of the next day based on the optimization of the market utility (objective) function $U$ with optimization (decision) variables $z$ taking values in $R^{|z|}$, (where $|\bullet|$ denotes a cardinality of a set). Variables $z$ are decomposed into hourly nodal power injection/withdrawal variables $q$ and variables $u$, describing the hourly power flow in the system (for example, in DC model variables $u$ are identified with phase angles of nodal voltage levels, in full AC model $u$ include hourly power flow volumes in each line of the network, nodal voltage magnitudes and phase angles, etc.). For simplicity we assume that there is only one generation or consumption unit in each node and denote the corresponding nodal power injection/withdrawal volume in hour $h$ by $q^{n,h}$, $n=1,...,N$, with $N$ being total number of nodes in the system. Thus, $z=\{q,u\}$ and $q=\{q^{n,h}\}$. The feasible set of $z$ produced by constraints $\{G(z)\}$, involving both constraints in the form of equalities as well as those in the form of weak inequalities, is assumed to be nonempty compact subset of Euclidean space $R^{|z|}$. (To ensure compactness of the feasible set $z$ and avoid dealing with multiple physically equivalent solutions to DAM problem differing in phase angles values by $2\pi$ multiples, the set $\{G(z)\}$ includes constraints on phase angles values limiting their feasible values to the corresponding close intervals with widths less than $2\pi$).

The DAM optimization problem has the form
$$\max_{\substack{z, \\ s.t.\{G\}}} U(q) \qquad (2)$$

with market utility function $U = U^{cons} - U^{prod}$, where $U^{cons}$ and $U^{prod}$ are additively separable total daily cost of power consumption as bid by the consumers and the total daily cost of power as offered by the producers respectively. It is assumed that market players submit distinct bids/offers with respect to each consumption/generation unit. If consumers may submit only totally inelastic DAM bids, then usually the term $U^{cons}$ is omitted from the function $U$, and the set $\{G\}$ is properly extended to account for the fixed consumption volumes.

We assume that the constraint functions $\{G(z)\}$ are twice continuously differentiable functions of $z$ and the bid/offered prices are piece-wise continuous functions of the respective injection/withdrawal volumes, hence, the function $U(q)$ is continuous. Compactness of the feasible set and continuity of the objective function ensure that its maximum value is attainable. In what follows we also assume that (2) has a unique maximizer, which we denote by $z^*$.

The set $\{G\}$ typically includes transmission constraints due to the power flow thermal or security limits, nodal power balance equations, power losses and power flow equations (Kirchhoff laws), generating unit constraints. The latter account for minimal/maximal output volumes, ramping rates, fuel constraints, etc. Let's denote a subset of $\{G(q,u)\}$, binding at $z=z^*$, by $\{C(q,u)\}$ and partition the latter as follows
$$\{C(q,u)\} = \{C_q(q), C_{qu}(q,u), C_u(u)\} \qquad (3)$$

with constraint subset $\{C_q(q)\}$ describing generation and consumption unit constraints, which are assumed to be functions of injection/consumption volumes in that node only (possibly in different hours), subset $\{C_{qu}(q,u)\}$ being a set of hourly nodal power balance constraints of the form $C_{qu}^{n,h}(q,u) = f_{qu}^{n,h}(u) \pm q^{n,h} = 0$, with "plus" sign if $q^{n,h}$ corresponds to power consumption and with "minus" sign in the case of power generation, $|C_{qu}(q,u)| = 24N$, subset $\{C_u(u)\}$ being attributed to electric power flow equations, network power flow limits, etc., with subscripts "$q$", "$qu$", and "$u$" referring to the constraint type and not being vector or matrix indices. We assume that a set of binding constraints $\{C(q,u)\}$ has cardinality less than $|z|$ and satisfies constraint qualification condition with respect to variables $(q,u)$. In this setting the nodal price in node $n$ at hour $h$, which we denote by $\lambda_{n,h}$, equals the Lagrange multiplier to the constraint $C_{qu}^{n,h}(q,u) = 0$. We note that for different hours of the day the sets of binding inequality constraints can differ.

As the firm $g$ operates generating units located in a set of nodes $\{i\}$, let $\{s\}$ denote the rest of the nodes in the system, $\{x^{i,h}\}$ and $\{y^{s,h}\}$ denote $\{q^{i,h}\}$ and $\{q^{s,h}\}$ respectively, $x = \{x^{i,h}\}$, $y = \{y^{s,h}\}$, $U_g(x)$ denote the offered costs by $g$ and $U_{\bar{g}}(y)$ denote the other market players contribution to the market utility function, so that $U = U_{\bar{g}}(y) - U_g(x)$. The sets of constraints $\{C_q(q)\}$, $\{C_{qu}(q,u)\}$ assume the following structure:

$$\{C_q(q)\} = \{C_x(x), C_y(y)\}$$

$$\{C_{qu}(q,u)\} = \{C_{xu}(x,u), C_{yu}(y,u)\},$$

$C_{xu}^{i,h}(x,u) = f_{xu}^{i,h}(u) - x^{i,h} = 0$, $C_{yu}^{s,h}(y,u) = f_{yu}^{s,h}(u) \pm y^{s,h} = 0$.

The constraint set $\{\overline{C}\} \equiv \{C_y(y), C_{xu}(x,u), C_{yu}(y,u), C_u(u)\}$ is referred to as reduced binding constraint set at $z = z^*$. We observe that $x$ as well as $y$ (up to a sign) enter additively in the reduced constraint set functions.

Let's denote by $D_x$ the feasible set for $x$, i.e., the set of $x$ such that there exists $(y,u)$, so that $(x,y,u)$ satisfies the constraint set $\{G\}$, and by $D_{(y,u)}(x)$ the set of all those $(y,u)$ for a given value of $x \in D_x$. Fixing a feasible value of $x$ we obtain an optimization subproblem of (2):

$$W_{\overline{g}}(x) = \max_{\substack{y,u, \\ (y,u) \in D_{(y,u)}(x)}} U_{\overline{g}}(y) \quad (4)$$

with $W_{\overline{g}}(x)$ being the optimal value of the market utility function when the firm $g$ supplies $x$ costless volumes of power in the corresponding nodes at respective hours of the day. The reduced binding constraint set $\{\overline{C}\}$ specifies a set of binding constraints for the problem (4). We emphasize that the set $\{\overline{C}\}$ may not satisfy constraint qualification for a given solution for $(y,u)$ of subproblem (4), and can be potentially redundant with some constraints being functions of the others. Relation of subproblem (4) to the problem (2) is given by

$$\max_{\substack{z, \\ s.t.\{G\}}} U(q) = \max_{\substack{x, \\ x \in D_x}} [W_{\overline{g}}(x) - U_g(x)] \quad . \quad (5)$$

Note that some of the constraints $\{G\}$, like $\{C_x(x)\}$, impose restrictions on $x$ only and involve neither $y$ nor $u$. However, in general case the constraints $\{G\}$ may also imply some additional implicit constraints on $x$, including the binding ones. Hence, generally speaking, the set of binding constraints on variables $x$ only, implied by $\{C(q,u)\}$, may be larger than $\{C_x(x)\}$.

Let's introduce costless power injections $\{a^{i,h}\}$ with $a^{i,h}$ being injection in node $i$ at hour $h$ of the day. That formally corresponds to a substitution $x^{i,h} \to x^{i,h} + a^{i,h}$ in the nodal power balance constraint for a node $i$ at hour $h$. The substitutions $x^{i,h} \to x^{i,h} + a^{i,h}$ in $C_{xu}^{i,h}(x,u)$ generate $D_{(y,u)}(x) \to D_{(y,u)}(x+a)$ and produce $a$-dependent feasible set of $x$: $D_x \to D_x(a)$. Therefore, when $\{a^{i,h}\}$ are introduced, the problem (4) is considered for $x \in D_x(a)$ and the corresponding value function equals $W_{\overline{g}}(x+a)$. Hence, if the Envelope Theorem (as stated in Appendix) is applicable for (2) at $a=0$, then we have the following identity for the nodal price at node $i$ at hour $h$:

$$\lambda_{i,h} = \left( \frac{\partial}{\partial a^{i,h}} \max_{\substack{x, \\ x \in D_x(a)}} \left[ W_{\overline{g}}(x+a) - U_g(x) \right] \right)\bigg|_{a=0}, \quad (6)$$

with

$$W_{\overline{g}}(x+a) = \max_{\substack{y,u, \\ (y,u) \in D_{(y,u)}(x+a)}} U_{\overline{g}}(y). \quad (7)$$

Equations (6)-(7) state that in general case both $W_{\overline{g}}(x+a)$ and the binding constraints on $x$, corresponding to the feasible set $D_x(a)$, are dependent on $a$ with two of these circumstances leading to nodal price $\lambda_{i,h}$ dependence on bids/offers of the market players other than $g$, while the latter circumstance being responsible for $\lambda_{i,h}$ dependence on the firm $g$ offers.

Dependence of $\lambda_{i,h}$ on the bids/offers of market players needs some clarification. Deformation of a bid/offer of a market player may result in the corresponding change in DAM outcomes including the player's nodal consumption/injection volumes. Let's consider a set of all bids/offers of a market player which result in the same values of nodal injection/consumption volumes as outcomes of DAM. We will say that the nodal price $\lambda_{i,h}$ is independent from the market player bid/offer, if any bid/offer from the abovementioned set results in the same value of $\lambda_{i,h}$. In that sense, the function $W_{\overline{g}}(x)$ is independent from the firm $g$ offers from the abovementioned set as it directly follows from (4). We would like to state conditions for a set of nodal prices $\{\lambda_{i,h}\}$ in nodes $\{i\}$ for all hours of the day to be independent from the firm $g$ offers.

**Proposition 1**. *If the following assumptions, additional to explicit assumptions stated in this section, are met: there exists $x^* \in D_x$ such that*

1. *for any $x$ in $\varepsilon(x^*)$ - some neighborhood of $x = x^*$ - the problem (4) has a unique solution $y = y^*(x)$, $u = u^*(x)$, and these functions are continuous in $\varepsilon(x^*)$;*

2. *in some neighborhood of $y = y^*(x)$ the function $U_{\overline{g}}(y)$ is twice continuously differentiable;*

3. *the set of constraints of the problem (4), binding at the solution of (4) for any $x$ in $\varepsilon(x^*)$, is given by the reduced constraint set $\{\overline{C}\}$ with cardinality no higher than $|y| + |u|$;*

4. *the set of binding constraints $\{\overline{C}\}$ satisfies constraint qualification at the solution of (4) for $x = x^*$: rank of $\left( \partial_y \overline{C}^m, \partial_u \overline{C}^m \right)\big|_{y=y^*(x^*), u=u^*(x^*)}$ is maximal;*

5. *if the kernel of $\left( \partial_y \overline{C}^m, \partial_u \overline{C}^m \right)\big|_{y=y^*(x^*), u=u^*(x^*)}$ is nontrivial, then we assume that the Hessian matrix with respect to $(y,u)$ of the Lagrangian*

$$\overline{L}(x,y,u) = U_{\overline{g}}(y) - \sum_{\overline{m}} \overline{\Lambda}_{\overline{m}}(x)\overline{C}^{\,\overline{m}}(x,y,u) \quad \text{at} \quad x = x^*,$$

$y = y^*(x^*)$, $u = u^*(x^*)$, *restricted to the kernel, is invertible, where* $\{\overline{\Lambda}_{\overline{m}}(x)\}$ *denote Lagrange multipliers associated with the reduced set of binding constraints* $\{\overline{C}\}$ *(for* $\overline{m}$ *corresponding to* $C^{i,h}_{xu}(x,u)$ *we denote* $\overline{\Lambda}_{\overline{m}}(x)$ *by* $\overline{\lambda}_{i,h}(x)$ *); then there exists a neighborhood of* $x = x^*$ *such that for any* $x$ *in that neighborhood the following statements are true:*

a. *the hourly nodal prices* $\lambda_{i,h}$ *are independent from the firm $g$ offers;*

b. $\overline{\lambda}_{i,h}$ *are continuously differentiable functions of $x$ and matrix* $\partial \overline{\lambda}_{i,h}(x)/\partial x^{i',h'}$ *is symmetric under* $(i,h) \leftrightarrow (i',h')$;

c. *if in addition in some open neighborhood of* $(y^*(x^*), u^*(x^*))$ *the function* $U_{\overline{g}}(y)$ *is (weakly) concave and the set* $\{\overline{C}\}$ *is convex (viewed as functions of variables $y$ and $u$), then there is an open neighborhood of* $x = x^*$ *such that the matrix* $\partial \overline{\lambda}_{i,h}(x)/\partial x^{i',h'}$ *is negative (semi-)definite in that neighborhood.*

Proof: Constraint qualification and smoothness of constraint and objective functions in (4) allow to use Karush–Kuhn–Tucker method and introduce Lagrange multiplier for each constraint of the set $\{\overline{C}\}$. Since the set $\{\overline{C}\}$ satisfies constraint qualification at $x = x^*$, so it does in some $\varepsilon(x^*)$ [1] as well. Therefore, the set $\{\overline{C}\}$ doesn't imply any constraints on $x$ in $\varepsilon(x^*)$ and the problem (4) can be considered for unconstrained $x \in \varepsilon(x^*)$. Considerations similar to the ones made in Appendix imply that $\forall x \in \varepsilon(x^*)$ the first-order condition for (4) allows to uniquely determine both $y = y^*(x)$ and $u = u^*(x)$ as well as Lagrange multipliers $\{\overline{\Lambda}_{\overline{m}}(x)\}$, with the corresponding functions being continuously differentiable functions in $\varepsilon(x^*)$ independent from the firm $g$ offers. Hence, $\overline{\lambda}_{i,h}(x)$ are continuously differentiable functions in $\varepsilon(x^*)$. From analysis presented in Appendix it follows that the Envelope Theorem is applicable to the subproblem (4) and in $\varepsilon(x^*)$ we have

$$\overline{\lambda}_{i,h}(x) = \partial W_{\overline{g}}(x)/\partial x^{i,h}. \quad (8)$$

---

[1] As our analysis requires frequent referrals to an open neighborhood of $x = x^*$, which is some open subset of $\varepsilon(x^*)$ defined above, for simplicity we will refer to that subset as $\varepsilon(x^*)$ as well and will tacitly assume that $\varepsilon(x^*)$ has been properly redefined. All the neighborhoods considered in this paper are assumed to be open neighborhoods.

Thus, the set of nodal prices $\overline{\lambda}_{i,h}(x)$ can be viewed as components of the gradient of function $W_{\overline{g}}(x)$. Since the first-order conditions for (4) make subset of the first-order conditions for (2), we have $\lambda_{i,h} = \overline{\lambda}_{i,h}(x^*)$. Hence, the nodal prices $\{\lambda_{i,h}\}$ in nodes $\{i\}$ for all hours of the day are independent from firm $g$ offers (provided that the offers result in the same supply volumes $x$ as defined above). From (8) it follows that $\partial \overline{\lambda}_{i,h}(x)/\partial x^{i',h'} = \partial \overline{\lambda}_{i',h'}(x)/\partial x^{i,h}$ in $\varepsilon(x^*)$. If $U_{\overline{g}}(y)$ is weakly concave and the residual set of binding constraints is convex, then negative semi-definiteness of $\partial \overline{\lambda}_{i,h}(x)/\partial x^{i',h'}$ follows from considerations in Appendix as Hessian matrices of $U_{\overline{g}}(y)$ and $-\sum_{\overline{m}} \overline{\Lambda}_{\overline{m}}(x)\overline{C}^{\,\overline{m}}(x,y,u)$ with respect to variables $(y,u)$ are negative semi-definite (with equality constraints not contributing to the latter). Moreover, if $U_{\overline{g}}(y)$ is strictly concave, then $\partial \overline{\lambda}_{i,h}(x)/\partial x^{i',h'}$ is negative definite. Thus, the Proposition 1 is proved.

Assumptions 3 and 4 of Lemma 1 imply that none of $g$'s offers is marginal, which limits applicability of the lemma, while the other Lemma 1 assumptions are often prove to be valid in practical applications.

We also note that in the case of DC current model with (weakly) concave quadratic objective function, $\partial \overline{\lambda}_{i,h}(x)/\partial x^{i',h'}$ is constant matrix in $\varepsilon(x^*)$.

Although $x = x^*$ is feasible, we do not require $\varepsilon(x^*) \in D_x$ as there may not be such an open neighborhood. Thus, the problem (4) is formally considered for all points $x \in \varepsilon(x^*)$, even the unfeasible ones. In this sense $\{x^{i,h}\}$ - the corresponding components of $x$ - are unconstrained and independent as long as the points belong to $\varepsilon(x^*)$ and hence the partial derivatives over $x^{i,h}$ can be considered.

We also note, that if the Envelope Theorem is applicable to the original problem (2) as well, the introduction of $\{a^{i,h}\}$ doesn't alter the set of binding constraints $\{C\}$: in some neighborhood of $a = 0$ no new binding constraints are introduced and no binding constraint becomes nonbinding. Since the reduced binding constraint set $\{\overline{C}\}$ doesn't impose any additional constraints on variables $x$, at $x = x^*(a)$ the set $\{C_x(x)\}$ is the set of binding constraints corresponding to the feasible set $D_x(a)$ and satisfying constraint qualification as a subset of $\{C(q,u)\}$. For introduction of $\{a^{i,h}\}$ doesn't introduce explicit dependence on $a$ in $\{C_x(x)\}$, using (6) and the statement of Lemma 2 in Appendix, we have

$$\lambda_{i,h} = \left( \frac{\partial}{\partial a^{i,h}} W_{\overline{g}}(x+a) \right)\bigg|_{\substack{x=x^*(0) \\ a=0}} = \left( \frac{\partial}{\partial x^{i,h}} W_{\overline{g}}(x) \right)\bigg|_{x=x^*} , \quad (9)$$

where we have used $x^* = x^*(0)$ and the fact that $x^*(a)$ - the solution for $x$ as a function of $a$ - is unique (which is implied by the assumption of the Envelope Theorem applicability for (2)).

Thus, if the problem (4) is considered only for $x \in D_x$, it may not be legitimate to take partial derivatives over $x^{i,h}$, thus derivatives are taken over $a^{i,h}$. However, due to the constraint qualification of $\{\overline{C}\}$, the problem (4) can be considered $\forall x \in \varepsilon(x^*)$, thus enabling computation of partial derivatives over $x^{i,h}$, which is manifest in (8) and (9). That observation also entails that $\overline{\lambda}_{i,h}(x)$ with $x = \{x^{i',h'}\}$ is interpreted as hourly nodal price in given hour $h$ and node $i \in \{i\}$, when at each hour $h'$ of the day and each node $i' \in \{i\}$ the firm $g$ injects costless power in the amount of $x^{i',h'}$.

## III. Generalization

In the previous section we dealt with the specific properties of the nodal prices in nodes $\{i\}$, containing generating units of the firm $g$, provided that conditions of the Proposition 1 hold for all hours of the day. However, as the network power flow can change during the day, the power flow constraints may imply some binding constraints on variables $x$ only in a few hours of the day. To address that issue, we generalize the statement of the Proposition 1 to a subset of hours of the day to account for hour specific conditions.

Let's partition $I = \{i\} \times \{1;2;...;24\}$ - the set of possible values of composite index $(i,h)$ - into $I = \hat{I} \cup \tilde{I}$ with $\hat{I} \cap \tilde{I} = \emptyset$. That induces a split of the variables $x$ into $x = (\hat{x}, \tilde{x})$ with $\hat{x} = x^{i,h}$, $(i,h) \in \hat{I}$, $\tilde{x} = x^{i,h}$, $(i,h) \in \tilde{I}$. We note that for a particular node $i$ the variable $x^{i,h}$ can belong to $\hat{x}$ at a certain hour $h$ of a day and belong to $\tilde{x}$ at another hour of the day. Let's further assume that the constraints $\{C_x(x)\}$ split into $\{C_x(x)\} = \{C_{\hat{x}}(\hat{x}), C_{\tilde{x}}(\tilde{x})\}$. Analogously to the treatment in the previous section, we define by $D_{\tilde{x}}$ the feasible set for $\tilde{x}$, i.e., the set of $\tilde{x}$ such that there exists $(\hat{x}, y, u)$, so that $(x, y, u)$ satisfies the constraint set $\{G\}$, and by $D_{(\hat{x},y,u)}(\tilde{x})$ the set of all those $(\hat{x}, y, u)$ for a given value of $\tilde{x} \in D_{\tilde{x}}$. Let's define $\hat{U}_g(\hat{x})$ and $\tilde{U}_g(\tilde{x})$ - firm $g$ offer costs function for $\hat{x}$ and $\tilde{x}$ respectively – through $U_g(x) = \hat{U}_g(\hat{x}) + \tilde{U}_g(\tilde{x})$, and define $\{\overline{C}\}$ as $\{\overline{C}\} \equiv \{C_{\hat{x}}(\hat{x}), C_y(y), C_{xu}(x,u), C_{yu}(y,u), C_u(u)\}$. We also define:

$$W_{\overline{g}}(\tilde{x}) = \max_{\substack{\hat{x},y,u, \\ (\hat{x},y,u) \in D_{(\hat{x},y,u)}(\tilde{x})}} [U_{\overline{g}}(y) - \hat{U}_g(\hat{x})] \quad (10)$$

Introducing costless power injections $\{a^{i,h}\}$ for $(i,h) \in \tilde{I}$, we arrive at the following generalization of the Proposition 1.

**Proposition 2**. *If the following assumptions, additional to explicit assumptions stated in section II, are met: there exists $\tilde{x}^* \in D_{\tilde{x}}$ such that*

*1. for any $\tilde{x}$ in $\varepsilon(\tilde{x}^*)$ - some neighborhood of $\tilde{x} = \tilde{x}^*$ - the problem (10) has a unique solution $\hat{x} = \hat{x}^*(\tilde{x})$, $y = y^*(\tilde{x})$, $u = u^*(\tilde{x})$, and these functions are continuous in $\varepsilon(\tilde{x}^*)$;*

*2. in some neighborhood of $\hat{x} = \hat{x}^*(\tilde{x})$, $y = y^*(\tilde{x})$ the function $U_{\overline{g}}(y) - \hat{U}_g(\hat{x})$ is twice continuously differentiable;*

*3. the set of constraints of a problem (10), binding at the solution of (10) for any $\tilde{x}$ in $\varepsilon(\tilde{x}^*)$, is given by the reduced constraint set $\{\overline{C}\}$ with cardinality no higher than $|\hat{x}| + |y| + |u|$;*

*4. the set of binding constraints $\{\overline{C}\}$ satisfies constraint qualification at the solution of (4) for $\tilde{x} = \tilde{x}^*$: rank of $\left(\partial_{\hat{x}} \overline{C}^{\overline{m}}, \partial_y \overline{C}^{\overline{m}}, \partial_u \overline{C}^{\overline{m}}\right)\big|_{\hat{x}=\hat{x}^*(\tilde{x}^*), y=y^*(\tilde{x}^*), u=u^*(\tilde{x}^*)}$ is maximal;*

*5. if the kernel of $\left(\partial_{\hat{x}} \overline{C}^{\overline{m}}, \partial_y \overline{C}^{\overline{m}}, \partial_u \overline{C}^{\overline{m}}\right)\big|_{\hat{x}=\hat{x}^*(\tilde{x}^*), y=y^*(\tilde{x}^*), u=u^*(\tilde{x}^*)}$ is nontrivial, then we assume that the Hessian matrix with respect to $(\hat{x}, y, u)$ of the Lagrangian $\overline{L}(x,y,u) = U_{\overline{g}}(y) - \hat{U}_g(\hat{x}) - \sum_{\overline{m}} \overline{\Lambda}_{\overline{m}}(\tilde{x}) \overline{C}^{\overline{m}}(x,y,u)$ at $\tilde{x} = \tilde{x}^*, \hat{x} = \hat{x}^*(\tilde{x}^*), y = y^*(\tilde{x}^*), u = u^*(\tilde{x}^*)$, restricted to the kernel is invertible, where $\{\overline{\Lambda}_{\overline{m}}(\tilde{x})\}$ denote Lagrange multipliers associated with the reduced set of binding constraints $\{\overline{C}\}$ (for $\overline{m}$ corresponding to $C_{\tilde{x}u}^{i,h}(\tilde{x},u)$ we denote $\overline{\Lambda}_{\overline{m}}(\tilde{x})$ by $\overline{\lambda}_{i,h}(\tilde{x})$);*

*then there exists an open neighborhood of $\tilde{x} = \tilde{x}^*$ such that for any $\tilde{x}$ in that neighborhood the following statements are true:*

*a. the hourly nodal prices $\overline{\lambda}_{i,h}$ for $(i,h) \in \tilde{I}$ are independent from the firm $g$ offers in $(i,h) \in \tilde{I}$;*

*b. $\overline{\lambda}_{i,h}$ for $(i,h) \in \tilde{I}$ are continuously differentiable functions of $\tilde{x}$ and matrix $\partial \overline{\lambda}_{i,h}(\tilde{x}) / \partial \tilde{x}^{i',h'}$ is symmetric under $(i,h) \leftrightarrow (i',h')$ for $(i,h),(i',h') \in \tilde{I}$;*

*c. if in addition in some open neighborhood of $\hat{x} = \hat{x}^*(\tilde{x}^*)$, $y = y^*(\tilde{x}^*)$, $u = u^*(\tilde{x}^*)$ the function $U_{\overline{g}}(y) - \hat{U}_g(\hat{x})$ is (weakly) concave and the set $\{\overline{C}\}$ is convex (viewed as functions of variables $(\hat{x}, y, u)$), then there is an open neighborhood of $\tilde{x} = \tilde{x}^*$ such that the matrix $\partial \overline{\lambda}_{i,h}(\tilde{x}) / \partial \tilde{x}^{i',h'}$ for $(i,h),(i',h') \in \tilde{I}$ is negative (semi-)definite in that neighborhood.*

The proof of the Proposition 2 is fully analogous to the one of Proposition 1.

## IV. EXAMPLE

In this section we consider an example illustrating the findings in the case of a model involving binding intertemporal constraint. Consider one node two generating units model with DAM solved for a time period composed of two consecutive hours $h = \{1;2\}$. Let demand be given by $\alpha(1.5n - y_1^D)$ and $\alpha(3n - y_2^D)$ in the first and second hour respectively with positive parameter $\alpha$ having units of $\$/(MWh)^2$, while nonnegative parameter $n$ and optimization variables $y_1^D$, $y_2^D$ having units of power volume. Let the firm $g$ supply $x_h$ units of power in hour $h$, and the other firm $g'$, operating a generator with maximal output $4n$ and upward ramping constrained by $n$, offer a price $7\alpha n/4$ to supply volume of power $y_1^S$ in the first hour, and a price $\alpha n$ to supply volume of power $y_2^S$ in the second hour. The utility function for a subproblem (4) is given by

$$U_{\bar{g}}(y) = \alpha[1.5ny_1^D - (y_1^D)^2/2] + \alpha[3ny_2^D - (y_2^D)^2/2] - 7\alpha n y_1^S/4 - \alpha n y_2^S$$

with the following constraints:

$$0 \leq y_1^D \leq 1.5n, \ 0 \leq y_2^D \leq 3n, \ 0 \leq y_h^S \leq 5n, \ x_h + y_h^S = y_h^D, \ y_2^S \leq y_1^S + n.$$

For sufficiently small positive $x_h^*$, $h = \{1;2\}$, there is an open neighborhood $\varepsilon(x^*)$ such that feasible set for variables $\{y_1^D, y_2^D, y_1^S, y_2^S\}$ is nonempty, the maximizer is unique and continuous in $(x_1, x_2)$:

$$y_1^D(x_1, x_2) = 3n/8 + (x_1 - x_2)/2, \ y_2^D(x_1, x_2) = 11n/8 + (x_2 - x_1)/2,$$

$$y_1^S(x_1, x_2) = 3n/8 - (x_1 + x_2)/2, \ y_2^S(x_1, x_2) = 11n/8 - (x_1 + x_2)/2,$$

with hourly prices

$$\bar{\lambda}_1(x_1, x_2) = \alpha[9n/8 + (x_2 - x_1)/2], \ \bar{\lambda}_2(x_1, x_2) = \alpha[13n/8 + (x_1 - x_2)/2]. \quad (11)$$

We note that due to the presence of the binding intertemporal constraint the hourly nodal prices are not given by the intersection of supply and demand curves in that hour. Also, in the first hour the nodal price is below the average variable cost of the firm $g'$ and it suffers a financial loss in that hour. However, in the second hour the firm $g'$ operates profitably with high margin. Due to the ramping constraint it is overall profitable for the firm $g'$ to produce some power in the first hour (incurring financial losses) to be able to produce more power in the second hour (with a profit exceeding the losses of the previous hour). As $g$ increases injection in the first hour, it lowers the nodal price in that hour and makes $g'$ lower output in the first hour, which – due to the binding ramping constraint – results in lower production by $g'$ in the second hour as well, thus increasing the nodal price in the second hour to compensate for higher economic losses per unit of power in the first hour. That is the mechanism of the nodal price increase in the second hour due to increase of power injection by $g$ in the first hour.

Mathematically, in $\varepsilon(x^*)$ the upward ramping constraint is the only binding inequality constraint, and the set of binding constraints $\bar{C} = \{x_h + y_h^S = y_h^D, y_2^S = y_1^S + n\}$ is unaltered and satisfies the constraint qualification at the solution. The kernel of $(\partial_y \bar{C})|_{y=y^*(x)}$ is spanned by $\delta y_1^D = \delta y_2^D = \delta y_1^S = \delta y_2^S$, and hence the Hessian of the Lagrangian function is invertible in the tangent space of the surface, defined by the binding constraints. Therefore, assumptions of the Proposition 1 hold and we have

$$\partial \bar{\lambda}_1 / \partial x_2 = \partial \bar{\lambda}_2 / \partial x_1,$$

which agrees with the direct computation using (11). This examples shows that nodal price in the first hour depends on the firm $g$ injection in the second hour and vise versa. Thus, consideration of binding intertemporal constraints is important in multi-period DAM model. Also, the 2x2 matrix nodal price response matrix is negative semi-definite with eigenvalues $0, (-\alpha)$, the null vector and the eigenvector are given by $(1,1)$ and $(1,-1)$ respectively. In the given example changing production volumes by the same amount in each hour will not change the nodal prices, whereas increasing the injection in one hour and decreasing it by the same amount in the other hour is the most efficient way to exploit the sensitivity of nodal prices with respect to the firm $g$ production volumes (that also directly follows from (11) as the firm's injections enter the nodal prices only in $(x_2 - x_1)$ combination).

## V. CONCLUSION

In this paper we studied properties of the nodal price response in DAM with respect to infinitesimal injections of costless power by a firm $g$, operating generating units assigned to different nodes. The constraint set of the DAM optimization problem includes (possibly, intertemporal) constraints on power consumption/production volumes specified by the market players and constraint set originating from the network (power flow equations, limits on power flows, etc.). We have shown that if the reduced set of binding constraint (defined as the full set of binding constraints excluding $g$'s generating units constraints such as minimal/maximal output, ramping, etc.) doesn't induce any additional binding constraints on the production volumes by the firm $g$, then (given validity of the other assumptions, stated in Proposition 1) the nodal price response matrix is symmetric and negative (semi-)definite for all hours of the day: matrix $\partial \bar{\lambda}_{i,h}(x)/\partial x^{i',h'}$ is symmetric under $(i,h) \leftrightarrow (i',h')$. Symmetry property of the nodal price response matrix is analogous to the statement of the reciprocity theorem for networks, which implies that if electromagnetic force insertion in one loop produces a current in another loop, then insertion by the electromagnetic force in the position of the current produces equal current in the first loop. That result is further extended in Proposition 2 for the case, when additional binding constraints are absent only for a subset of hours of the day. These findings can be viewed as generalization of results obtained in [5], [9] to the case of power systems with nonlinear and intertemporal constraints.

## VI. Appendix

In this section we recall some basic properties of the Lagrange multiplier derivatives over exogenous parameters in constrained optimization problem

$$\max_{\substack{v,\\ s.t.\{G(x,v)\}}} F(v) \quad (12)$$

with $x \in X$, $v \in \Omega$, where $X$ and $\Omega$ are open subsets of Euclidean spaces $R^{|x|}$ and $R^{|v|}$ respectively, $F(v)$ is twice continuously differentiable objective function defined on $\Omega$, constraints $\{G(x,v)\}$ are defined on $X \times \Omega$ and specify equality as well as weak inequality constraints. Functions $G(x,v)$ are assumed to be twice continuously differentiable functions of the variables $(x,v)$ with $x$ specifying a set of exogenous parameters $x = \{x^i\}$, $i=1,...,|x|$. The feasible set defined by $\{G(x,v)\}$ for each $x \in X$ is assumed to be nonempty compact (bounded and closed) set in $\Omega$. The extreme value theorem implies that for each $x$ there exists at least one solution for $v$ belonging to the feasible set and maximizing (12). We emphasize that under the stated assumptions even if the objective function and constraints are smooth functions and there is unique maximizer of (12), the value function and maximizer of (12) in general case need not be continuous functions of the parameters $x$.

We will make the following additional assumptions on (12): there exists $x^* \in X$ such that

a) for any $x$ in $\varepsilon(x^*)$ - an open neighborhood of $x^*$, $\varepsilon(x^*) \in X$, - the problem (12) has unique solution, denoted as $v = v^*(x)$, and function $v^*(x)$ is continuous in $\varepsilon(x^*)$;

b) the set of binding constraints is unaltered in $\varepsilon(x^*)$ (we denote by $\{C^m(x,v)\}$, $m=1,...,M$, the subset of $\{G(x,v)\}$ representing constraints, binding at $x = x^*$, $v = v^*(x^*)$);

c) $|x| = M \leq |v|$ and all functions $C^m(x,v)$ have the form $C^m(x,v) = c^m(v) - x^m$;

d) rank of $\partial_v C^m(x,v)\big|_{v=v^*(x^*)}$ is maximal;

e) if the kernel of $\partial_v C^m(x,v)\big|_{v=v^*(x^*)}$ is nontrivial, then we assume that the Hessian matrix with respect to $v$ of the Lagrangian function $L(x,v) = F(v) - \sum_m \Lambda_m(x) C^m(x,v)$ at $x = x^*$, $v = v^*(x^*)$, i.e., $D_v^2 L\big|_{\substack{v=v^*(x^*)\\x=x^*}}$, restricted to the kernel of $\partial_v C^m(x,v)\big|_{v=v^*(x^*)}$, is invertible (where $\{\Lambda_m(x)\}$ denote Lagrange multipliers associated with constraints $\{C^m(x,v)\}$).

A few consequences readily follow: assumption on twice continuous differentiability of $G(x,v)$ implies that functions $c^m(v)$ are also twice continuously differentiable, which – given continuity of $v^*(x)$ in $\varepsilon(x^*)$ - implies continuity of $\partial_v C^m(x,v)\big|_{v=v^*(x^*)}$ in $\varepsilon(x^*)$; the latter coupled with assumption "d" entails, that there is a neighborhood of $x^*$, such that rank of $\partial_v C^m(x,v)\big|_{v=v^*(x^*)}$ is maximal in that neighborhood. Likewise, there is an open neighborhood of $x^*$, such that assumption "e" is valid for all $x$ in that neighborhood

The Lagrange method is applicable in some $\varepsilon(x^*)$, (the footnote 1 is applicable to the present section as well), and we have the following set of equations, representing the first-order necessary condition for a critical point $\forall x \in \varepsilon(x^*)$:

$$\begin{cases} \partial_v L = 0 \\ -C^m(x,v) = 0, m=1,...,M \end{cases}, \quad (13)$$

supplemented with nonnegativity conditions for Lagrange multipliers associated with binding inequality constraints. Given the function $v^*(x)$, maximality of $\partial_v C^m(x,v)$ rank in $\varepsilon(x^*)$ allows to solve $\partial_v L = 0$ for $\Lambda_m = \Lambda_m^*(x)$, $m=1,...,M$, in some $\varepsilon(x^*)$. Continuity of $v^*(x)$ and twice continuous differentiability of $F(v)$ as well as $C^m(x,v)$, $m=1,...,M$, imply that $\Lambda_m^*(x)$ are continuous in that $\varepsilon(x^*)$.

It follows from [11] that rank of $\partial_v C^m(x,v)\big|_{v=v^*(x^*)}$ is maximal and $D_v^2 L\big|_{\substack{v=v^*(x^*)\\x=x^*}}$, restricted to the kernel of $\partial_v C^m(x,v)\big|_{v=v^*(x^*)}$, is invertible (with the latter condition present only if the kernel is nontrivial) if and only if the bordered Hessian matrix defined by

$$H = \begin{bmatrix} D_v^2 L & (\partial_v C)^T \\ (\partial_v C) & 0 \end{bmatrix}_{(x^*, v^*(x^*))} \quad (14)$$

is invertible. We note that since determinant of $H$ is nonzero at $x = x^*$ and is continuous function of $x$, it is also nonzero in some $\varepsilon(x^*)$.

It is well known, [12], that nondegeneracy preserving substitution $C^m(x,v) \to -C^m(x,v)$ in the bordered Hessian matrix produces Jacobian matrix for a set of equations (13). Hence implicit function theorem implies that in some neighborhood of $(x = x^*, v = v^*(x^*), \Lambda_m = \Lambda_m^*(x^*))$ there exists a unique (continuously differentiable with respect to $x$) solution of (13). We stress that in general case that doesn't imply that the solution to (12) is unique and/or continuously differentiable as the implicit function theorem provides a local solution to the necessary condition for extremum, while the optimization problem might have multiple and/or discontinuous solutions. However in our case assumptions of uniqueness and continuity of solution $v = v^*(x)$ of (12) in some neighborhood of $x = x^*$ implies that that solution is the one identified by the implicit function theorem, and there is $\varepsilon(x^*)$ such that $v^*(x)$ is continuously differentiable function in that neighborhood. Since $F(v)$ is twice continuously differentiable function and in the neighborhood $\varepsilon(x^*)$ the maximizer of (12) $v = v^*(x)$ is unique and continuously differentiable, the value function for (12) defined by $V(x) = F(v^*(x))$ is also continuously differentiable in $\varepsilon(x^*)$.

Note that we can remove the assumption "e", if the assumption "a" is extended to require that $v^*(x)$ is continuously differentiable functions in $\varepsilon(x^*)$. In that case, the continuous differentiability of $\Lambda_m^*(x)$ follows directly from $\partial_v L = 0$ utilizing constraint qualification, which entails that in some $\varepsilon(x^*)$ the matrix

$$\left(\sum_{s=1,...,|v|} \partial_{v^s} C^m(x,v) \partial_{v^s} C^{m'}(x,v)\right)\bigg|_{v=v^*(x)}$$ is invertible.

Constraint qualification and $M \leq |v|$ also imply that a set $\{C^m(x,v)\}$ doesn't not impose any constraints on the possible values of $x$ (even for a case when the constraints do not have the form $C^m(x,v) = c^m(v) - x^m$): there is no continuously differentiable function $S(C,x)$ in some neighborhood of $(x, v = v^*(x))$, such that it is a function of $x$ only: $\partial_v S(C,x) = 0$ with $d_x S(C,x) \not\equiv 0$ and $\sum_m (\partial S/\partial C^m)^2 \not\equiv 0$. That conclusion is compatible with the constant rank theorem.

The continuity of the function $v^*(x)$ in some open neighborhood of $x^*$, as stated in assumption "a", can be replaced by more fundamental assumptions on the constraints functions $\{G(x,v)\}$. One such example is given by the Berge's Maximal Theorem in the framework of constraints generated correspondence.

Economic interpretation of the Lagrange multiplies as shadow prices [13] is given by the Envelope Theorem [14], which we reformulate below in a way suitable for our purposes.

**The Envelope Theorem**: *If* $f(z,a)$, $g^1(z,a),...,g^P(z,a), h^1(z,a),...,h^T(z,a)$, *are continuously differentiable functions of* $(z,a)$ *on* $R^{|z|} \times R^{|a|}$ *and there is* $D_a$ - *an open subset of* $R^{|a|}$, *such that for any* $a \in D_a$

- $z^*(a)$ - *is the unique solution of the optimization problem*

$$V(a) = \max_{\substack{z \in R^{|z|}, \\ s.t. \\ g^p(z,a)=0, \\ h^t(z,a) \leq 0}} f(z,a) \quad (15)$$

*with* $p = 1,...,P$, $t = 1,...,T$;

- $z^*(a)$ *is continuously differentiable function of* $a$ *in* $D_a$;

- *in* $D_a$ *the set of constraints binding at* $z^*(a)$ *is unaltered, has cardinality no higher than* $|z|$, *and satisfies the constraint qualification condition. We denote the binding constraints by* $C^l(z,a)$ *and the corresponding Lagrange multipliers by* $\Lambda_l(a)$, $l = 1,...,L$, $P \leq L \leq P + T$;

*then for any* $a \in D_a$ *the value function* $V(a)$ *is continuously differentiable function of* $a$ *and satisfies*

$$\frac{\partial}{\partial a}V(a) = \left(\frac{\partial}{\partial a}f(z,a)\right)\bigg|_{z=z^*(a)} - \sum_{l=1,...,L}\Lambda_l^*(a)\left(\frac{\partial}{\partial a}C^l(z,a)\right)\bigg|_{z=z^*(a)}, \quad (16)$$

where $\Lambda_l^*(a)$ denote the solution for the Lagrange multiplier $\Lambda_l$, $l = 1,...,L$.

The proof of continuous differentiability of the value function readily follows from uniqueness of the maximizer and continuous differentiability of the objective function and the maximizer, while application of the chain rule to $\partial f(z^*(a),a)/\partial a$ together with the first-order necessary condition, stability of the binding constraint set in $D_a$, and utilization of $\partial C^l(z^*(a),a)/\partial a = 0$, $l = 1,...,L$, yields (16). We also note that uniqueness of the solution for $z^*(a)$ and the constraint qualification imply uniqueness of the solution for $\Lambda_l^*(a)$ in $D_a$.

In some formulations of the Envelope Theorem continuous differentiability of $\Lambda_l^*(a)$, associated with inequality constraints, is also required but condition of stability of the binding constraint set in $D_a$ is omitted.

We note that if the set of binding constraints is unaltered in $D_a$ and constraint qualification is satisfied at some point $a = a^* \in D_a$, then there is an open set $D_a' \subset D_a$, such that the constraint qualification holds for any $a \in D_a'$.

Let's denote by $D_z(a)$ the feasible set of $z$, defined by the constraints in (15), and partition variables $z$ into $z = (x,v)$, and let $D_x(a)$ be the feasible set for $x$, i.e., the set of all $x$ such that there exists $v$, such that $z = (x,v)$ belongs to $D_z(a)$, and let $D_v(x,a)$ be the set of all those $v$ for a given value of $x \in D_x(a)$. In the correspondence formalism $D_x(a)$ and $D_v(x,a)$ are respectively identified as domain and image of $x$ under (parameterized by $a$) correspondence, defined by map of feasible values of $x$ into a subset of $R^{|v|}$, which contains the union of $D_v(x,a)$ over all $x \in D_x(a)$.

Let's define

$$Q(x,a) = \max_{\substack{v \in R^{|v|}, \\ v \in D_v(x,a)}} f(z,a), \quad (17)$$

then we have

$$V(a) = \max_{\substack{x \in R^{|x|}, \\ x \in D_x(a)}} Q(x,a). \quad (18)$$

**Lemma 1.** *Let assumptions of The Envelope Theorem hold for optimization problem (15), the variables* $z$ *and the constraints, binding in* $D_a, C^l(z,a)$, $l = 1,...,L$, *can be* *partitioned into sets* $z = (x, v)$ *and* $\{C^l\} = (\{\overline{C}^{\bar{l}}\}, \{\widetilde{C}^{\tilde{l}}\})$, $\bar{l} = 1,...,\overline{L}$, $\tilde{l} = 1,...,\widetilde{L}$, *respectively, such that*

1. $\{\widetilde{C}^{\tilde{l}}\}$ *are independent from both* $v$ *and* $a$, $\widetilde{C}^{\tilde{l}} = \widetilde{C}^{\tilde{l}}(x)$;

2. $\overline{L} \leq |v|$ *and rank of* $\partial \overline{C}^{\bar{l}}(x,v,a)/\partial v\big|_{z=z^*(a)}$ *is maximal in* $D_a$;

3. *for any* $(x,a)$ *such that* $x \in D_x(a)$, $a \in D_a$, *the maximizer for optimization problem (17) is unique and continuously differentiable function with respect to* $(x,a)$;

*then*

$$\frac{\partial}{\partial a} V(a) = \left(\frac{\partial}{\partial a} Q(x,a)\right)\bigg|_{x=x^*(a)}, \forall a \in D_a.$$

Proof: From assumptions of the Lemma 1 it follows that the function $Q(x,a)$ is continuously differentiable with respect to $(x,a)$ for $x \in D_x(a)$, $a \in D_a$. Since maximizer for (15) is unique and continuously differentiable function of $a$, so it is for (18). As $\{\overline{C}^{\bar{l}}\}$ do not induce any constraints on $x$ due to assumption 2 of the Lemma 1, the constraints $\widetilde{C}^{\tilde{l}} = \widetilde{C}^{\tilde{l}}(x)$ are the only binding constraints corresponding to $D_x(a)$. Hence, the Envelope Theorem is applicable to (18), which yields the statement of the Lemma 1.

**Lemma 2.** *Let assumptions of The Envelope Theorem hold for optimization problem (15), objective function* $f$ *is independent from* $a$, *the variables* $z$ *and the constraints, binding in* $D_a$, $C^l(z,a)$, $l=1,...,L$, *can be partitioned into* $z=(x,v)$ *and* $\{C^l\} = (\{\overline{C}^{\bar{l}}\}, \{\widetilde{C}^{\tilde{l}}\})$, $\bar{l} = 1,...,\overline{L}$, $\tilde{l} = 1,...,\widetilde{L}$, *respectively, such that* $|a| = |x|$ *and*

1. *the objective function is additively separable:* $f = \bar{f}(v) + \tilde{f}(x)$;

2. $\{\widetilde{C}^{\tilde{l}}\}$ *are independent from both* $v$ *and* $a$, $\widetilde{C}^{\tilde{l}} = \widetilde{C}^{\tilde{l}}(x)$, $\tilde{l} = 1,...,\widetilde{L}$;

3. *variables* $x$ *and* $a$ *enter* $\{\overline{C}^{\bar{l}}\}$ *in such a way that* $\overline{C}^{\bar{l}} = \overline{C}^{\bar{l}}(x+a, v)$, $\bar{l} = 1,...,\overline{L}$;

4. $\overline{L} \leq |v|$ *and rank of* $\partial \overline{C}^{\bar{l}}(x,v,a)/\partial v\big|_{z=z^*(a)}$ *is maximal in* $D_a$;

5. *for any* $(x,a)$ *such that* $x \in D_x(a)$, $a \in D_a$, *the maximizer for optimization problem*

$$W(x,a) = \max_{\substack{v \in R^{|v|}, \\ v \in D_v(x,a)}} \bar{f}(v) \quad (19)$$

*is unique and continuously differentiable with respect to* $(x,a)$;

*then* $W = W(x+a)$ *and* $\partial V(a)/\partial a = (\partial W(x+a)/\partial x)\big|_{x=x^*(a)}$, $\forall a \in D_a$. *Also, if* $\{a = 0\} \in D_a$, *then* $\partial V(a)/\partial a\big|_{a=0} = (\partial W(x)/\partial x)\big|_{x=x^*(0)}$.

Proof: Both continuous differentiability of $\bar{f}(v)$ and maximizer for (19) as well as uniqueness of the latter ensure that $W(x,a)$ is continuously differentiable function of both variables. Lemma 2 assumptions imply that KKT method is applicable to the problem (19) and since $a$ and $x$ enter the first-order necessary conditions only through $(x+a)$, we have $W(x,a) = W(x+a)$. Application of the Lemma 1 completes the proof:

$$\frac{\partial}{\partial a} V(a) = \left(\frac{\partial}{\partial a} [W(x+a) + \tilde{f}(x)]\right)\bigg|_{x=x^*(a)} = \left(\frac{\partial}{\partial x} W(x+a)\right)\bigg|_{x=x^*(a)}, \forall a \in D_a,$$

with immediate result for the case $\{a = 0\} \in D_a$.

*Symmetry of* $\partial \Lambda_i^*(x)/\partial x^j$

Applicability of the Envelope Theorem for (12) follows from the considerations above and, hence,

$$\Lambda_i^*(x) = \partial F(v^*(x))/\partial x^i \text{ in } \varepsilon(x^*). \quad (20)$$

However, the function $\Lambda_i^*(x)$ is continuously differentiable in $\varepsilon(x^*)$, hence the value function $V(x) = F(v^*(x))$ is twice continuously differentiable in that neighborhood. That in turn implies

$$\partial \Lambda_i^*(x)/\partial x^j = \partial^2 F(v^*(x))/\partial x^i \partial x^j, \quad (21)$$

and, therefore, the square matrix $\partial \Lambda_i^*(x)/\partial x^j$ is symmetric in $\varepsilon(x^*)$. These results readily follows from $\sum_s \partial L/\partial v^s\big|_{\substack{v=v^*(x) \\ \Lambda = \Lambda^*(x)}} \partial v^{*s}(x)/\partial x^i = 0$ using

$$\frac{\partial C^m(x,v)}{\partial x^i}\bigg|_{v=v^*(x)} + \sum_s \frac{\partial C^m(x,v)}{\partial v^s}\bigg|_{v=v^*(x)} \frac{\partial v^{*s}(x)}{\partial x^i} = 0. \quad (22)$$

Symmetry of $\partial \Lambda_i^*(x)/\partial x^j$ can be also linked to the symmetry of Hessian matrix of $L$ with respect to variables $v$:

$$\frac{\partial \Lambda_i^*(x)}{\partial x^j} = \sum_{s,k} \frac{\partial^2 L}{\partial v^s \partial v^k}\bigg|_{\substack{v=v^*(x) \\ \Lambda = \Lambda^*(x)}} \frac{\partial v^{*s}(x)}{\partial x^i} \frac{\partial v^{*k}(x)}{\partial x^j} \quad (23)$$

in $\varepsilon(x^*)$, which implies that the square matrix $\partial \Lambda_i^*(x)/\partial x^j$ is symmetric. Eq. (23) also entails

$$\frac{\partial \Lambda_i^*(x)}{\partial x^j} = \sum_{s,k} \frac{\partial^2 F(v)}{\partial v^s \partial v^k}\bigg|_{v=v^*(x)} \frac{\partial v^{*s}(x)}{\partial x^i} \frac{\partial v^{*k}(x)}{\partial x^j} \\ - \sum_{m,s,k} \Lambda_m^*(x) \frac{\partial^2 C^m(x,v)}{\partial v^s \partial v^k}\bigg|_{v=v^*(x)} \frac{\partial v^{*s}(x)}{\partial x^i} \frac{\partial v^{*k}(x)}{\partial x^j} \quad (24)$$

*Sign definiteness of* $\partial \Lambda_i^*(x)/\partial x^j$

Since in $\varepsilon(x^*)$ the function $v = v^*(x)$ is solution to (12), and bordered Hessian is invertible, the square matrix $\partial^2 L/\partial v^s \partial v^k\big|_{v=v^*(x)}$ defines negative definite form in the kernel of $\partial_v C^m(x,v)$, i.e., for $\delta v$ preserving the constraints. However, as $\partial v^{*s}(x)/\partial x^i$ doesn't belong to the kernel due to the term $\partial C^m(x,v)/\partial x^i\big|_{v=v^*(x)}$ in equation (22), in general case (23) doesn't imply that $\partial \Lambda_i^*(x)/\partial x^j$ is negative definite or negative semi-definite matrix. We note for the case of convex optimization problem the negative semi-definiteness of $\partial \Lambda_i^*(x)/\partial x^j$ in $\varepsilon(x^*)$ follows from (24). Also, (22) and assumed structure of $\{C^m(x,v)\}$ entail that $\sum_i (\partial v^{*s}(x)/\partial x^i)\rho^i \neq 0$ with any non-zero vector $\rho$. Hence, for the cases of convex optimization problem with strictly concave function $F(v)$, the matrix $\partial \Lambda_i^*(x)/\partial x^j$ is negative definite in $\varepsilon(x^*)$.

## VII. Acknowledgment

## IX. Biography

**Vadim Borokhov**, received M.S. (1998) in Physics from Moscow Institute of Physics and Technology, and Ph.D. (2004) in Physics from California Institute of Technology. In 2004-2006 he worked as consultant at LLC "Carana" (Moscow, Russia) on deregulation of the Russian electricity market, in 2006-2008 worked as utilities equity analyst at Renaissance Capital (Moscow, Russia), in 2008-2014 served as a project manager and department head at the Market Council (Moscow, Russia). Currently holds position of Advisor on power market development at En+ Development.